\documentclass[MSNbibl,number,citesort,seceqn,dvips]{arxbj}
\usepackage{upgreek}

\aid{0}
\volume{19}
\issue{3}
\pubyear{2013}
\firstpage{886}
\lastpage{904}
\doi{10.3150/12-BEJ463} 

\makeatletter
\newcommand{\rrvert}{\vert}
\newcommand{\llvert}{\vert}
\newcommand{\COM}[1]{}

\newtheorem{theo}{Theorem}[section]
\newtheorem{lem}[theo]{Lemma}
\newcommand{\kb}[1]{\mathbf{#1}}
\newcommand{\vk}[1]{\kb{#1}}
\def\X{\vk{X}}
\newcommand{\ve}{\varepsilon}
\newcommand{\R}{\mathbb{R}}
\newcommand{\N}{\mathbb{N}}
\newcommand{\IF}{\infty}
\newcommand{\kal}[1]{{\mathcal{#1}}}
\newcommand{\sign}{\operatorname{sign}}
\def\X{\vk{X}}
\def\Y{\vk{Y}}
\def\SI{\Sigma}
\def\SIG{\SI}
\def\vks{\boldsymbol{\theta}}
\def\ZW#1{#1}
\def\U{\vk{U}}
\newcommand{\eqref}[1]{(\ref{#1})}
\newremark{remarks}{Remarks}
\newremark{example}{Example}
\newremark{remark}{Remark}
\newcommand{\fracf}[2]{({#1})/({#2})}
\newcommand{\fracc}[2]{{#1}/{(#2)}}
\makeatother

\begin{document}
\begin{frontmatter}

\title{Minima and maxima of elliptical arrays and spherical processes}
\runtitle{Extremes of elliptical arrays}

\begin{aug}
\author{\fnms{Enkelejd} \snm{Hashorva}\corref{}\ead[label=e1]{enkelejd.hashorva@unil.ch}}
\runauthor{E. Hashorva} 
\address{Faculty of Business and Economics, University of Lausanne, Extranef,
UNIL-Dorigny, 1015 Lausanne, Switzerland. \printead{e1}}
\end{aug}

\received{\smonth{6} \syear{2010}}
\revised{\smonth{5} \syear{2012}}

%
\begin{abstract}
In this paper, we investigate first the asymptotics of the minima of
elliptical triangular arrays. Motivated by the findings of
Kabluchko (\textit{Extremes} \textbf{14} (2011) 285--310), we discuss
further the asymptotic behaviour of the maxima
of elliptical triangular arrays with marginal distribution functions in
the Gumbel or Weibull max-domain of
attraction. We present an application concerning the asymptotics of the
maximum and the minimum of independent
spherical processes.
\end{abstract}

%
\begin{keyword}
\kwd{asymptotics of sample maxima}
\kwd{Brown--Resnick copula}
\kwd{Brown--Resnick process}
\kwd{Davis--Resnick tail property}
\kwd{Gaussian process}
\kwd{Penrose--Kabluchko process}
\kwd{spherical process}
\end{keyword}

\end{frontmatter}

\section{Introduction}\label{sec1}
It is well known that the maxima of Gaussian random vectors have
asymptotically independent components, a result going back to Sibuya
\cite{Sibuya1960}.
Recently, Kabluchko~\cite{Kabluchko2011} shows that the minima of the
absolute values of Gaussian random vectors have also asymptotically independent
components. The Gaussian framework is appealing from both theoretical
and applied point of view. In order to still consider Gaussian random
vectors for modelling asymptotically dependent risks, triangular arrays
of Gaussian random vectors with increasing dependence should be
considered -- this approach is suggested in H\"usler and Reiss \cite
{H1989}. As shown in the aforementioned paper, the maxima of Gaussian
triangular arrays
can be attracted by some max-stable distribution function with
dependent components which is referred to as the H\"usler--Reiss
distribution function. In fact,
the H\"usler--Reiss copula is a particular case of the
Brown--Resnick copula; a canonical example of a max-stable Brown--Resnick
process is first presented in Brown and Resnick~\cite{Brown1977}
in the context of the asymptotics of the maximum of Brownian motions.
See Kabluchko \textit{et al.}~\cite{Kabluchko2009} for the main properties of
Brown--Resnick
processes. Kabluchko~\cite{Kabluchko2011} discusses a more general
asymptotic framework analysing the maximum of independent Gaussian
processes showing that the
Brown--Resnick process
appears as the limit process if the underlying covariance functions
satisfy a certain asymptotic condition. Additionally,
the aforementioned paper investigates the asymptotics of the minimum of
the absolute value of independent Gaussian processes extending some previous
results of Penrose~\cite{Penrose1991}.

Indeed, Gaussian random vectors are a canonical example of elliptically
symmetric (for short elliptical) random vectors. Therefore, it is
natural to
consider Kabluchko's findings in the framework of elliptical random
vectors and spherical processes. Belonging to the class of
conditional Gaussian processes, spherical processes appear naturally in
diverse applications, see, for example, Falk \textit{et al.}~\cite{Falk2010},
or H\"usler \textit{et al.}~\cite{H2011a,H2011b}.

As shown in Hashorva~\cite{Hashorva2005,Hashorva2011} the maxima and
the minima (of absolute values) of elliptical random vectors
have asymptotically independent components. Elliptical random vectors
are defined by the marginal distribution functions and some
nonnegative definite matrix $\Sigma$, see \eqref{eq:el} below. If
$\SI_n,n\ge1$ are $k\times k$ correlation matrices pertaining to an elliptical
triangular array, the crucial condition
for the asymptotic behaviour of both maxima and minima is
%
\begin{eqnarray}
\label{eq:reis} \lim_{n \to\infty} c_n\bigl(\vk{1}
\vk{1}^\top- \Sigma_n\bigr) = \Gamma =:(
\gamma_{ij})_{i,j\le k} \qquad\mbox{with } \gamma_{ij}
\in(0,\IF), i\not=j, i,j\le k,
\end{eqnarray}
where $c_n,n\ge1$ is a sequence of positive constants determined by a
marginal distribution function of the elliptical random vectors, and
$\vk{1}=(1 ,\ldots,1)^\top\in\R^k$ (here $^\top$ stands for the
transpose sign).

In {Theorem~\ref{th:kabl:0}}, we specify the constants $c_n$ such that the
minima of absolute values of triangular arrays are attracted by some
min-infinitely divisible distribution function in~$\R^k$; the
dependence function of the limiting distribution function is indirectly
determined by the marginal distribution functions of the triangular
array. Utilising Kabluchko's approach, we reconsider the aforementioned
results for the maxima deriving some new representations for the
limiting distributions under the assumptions that the marginals of the
elliptical random vectors have distribution
function in the Gumbel or Weibull max-domain of attraction (MDA).

A direct application of our result concerns the asymptotics of maximum
and minimum (of absolute values) of independent spherical processes.
It turns out that the limiting process of the normalised maximum of
spherical processes is the same as that of Gaussian processes discussed
in Kabluchko~\cite{Kabluchko2011}, namely the max-stable Brown--Resnick
process. However, the norming constants are necessarily different. One
important consequence of our findings is that the Brown--Resnick process
is shown to be also the limit of the maximum of non-Gaussian processes.
When instead of maximum the minimum of absolute values of Gaussian
processes is considered, from the aforementioned reference, we know
that the limiting process is min-stable; we refer to that process as
Penrose--Kabluchko process. As demonstrated in our application,
Penrose--Kabluchko processes can be retrieved in the limit in the more
general framework of spherical processes.

The paper is organised as follows: Section~\ref{sec2} introduces our notation
and presents some preliminary results.
In Section~\ref{sec3}, we deal with the asymptotics of minima of absolute values
of elliptical triangular arrays. Section~\ref{sec4} investigates the maxima of
triangular arrays with marginal distribution functions in the MDA of
the Gumbel or the Weibull distribution. The applications mentioned above
are presented in Section~\ref{sec5}. Proofs of all the results are relegated to
Section~\ref{sec6}.

\section{Preliminaries}\label{sec2}
Let in the following $I,J$ be two non-empty disjoint index sets such
that $I\cup J=\{1 ,\ldots,k \},k\ge2,$
and define for $\vk{x}=(x_1 ,\ldots,x_k)^\top\in\R^k$ the
subvector of
$\vk{x}$ with respect to $I$ by
$\vk{x}_I=(x_i, i \in I)^\top$. If $\SI\in\R^{k\times k}$ is
a square matrix, then the matrix $\SI_{IJ}$ is obtained by retaining
both the rows and the columns of $\SI$ with indices in $I$ and in
$J$, respectively; similarly we define $\SI_{JI}, \SI_{JJ},
\SI_{II}$. Given $\vk{x},\vk{y}\in\R^k$ write
\begin{eqnarray*}
\vk{x}&>&\vk{y}\qquad \mbox{if } x_i>y_i,\forall i=1
,\ldots ,k,
\\
\vk{x}+ \vk{y}&=& (x_1+ y_1 ,\ldots,x_k+
y_k)^\top,\qquad c \vk{x}=(cx_1,
\ldots,cx_k)^\top,    c\in\R.
\end{eqnarray*}
The notation $\kal{B}_{a,b},a,b>0$ stands for a beta random variable
with probability density function
\[
\frac{\Gamma(a+b)}{\Gamma(a)\Gamma(b)}x^{a-1}(1- x)^{b-1}, \qquad x\in(0,1),
\]
where $\Gamma(\cdot)$ is the Euler Gamma function; $\vk{Y}\sim F$
means that the random vector $\vk{Y}$ has distribution function $F$.

Throughout this paper, $\vk{U}$ is a $k$-dimensional random vector
uniformly distributed on the unit sphere (with respect to the
$L_2$-norm) $\kal{S}_{k}$ of $\R^k$
being further independent of $R_k>0$ and $A,A_n,n\ge1$ are
$k$-dimensional square matrices such that $\SI=AA^\top$ and $\ZW{\SI_{n}=A_nA_n^\top}$
are positive definite correlation matrices (all entries in the main
diagonal are equal to 1). We write $\vk{U}_m$ if $m < k$ to mean again
that $\vk{U}_m$ has the uniform distribution on $\kal{S}_{m}$. The
distribution function of $R_k, k\ge1$ will be denoted by $H_k$,
whereas the distribution function of $R_kU_1$ will be denoted by $G$;
$\omega\in(0,\IF]$
is their common upper endpoint.

Let $\X=(X_1 ,\ldots,X_k)^\top, k\ge2$ be an elliptically symmetric
random vector with stochastic representation
%
\begin{eqnarray}
\label{eq:el} \X\stackrel{d} {=}R_k A \vk{U},
\end{eqnarray}
where $\stackrel{d}{=}$ stands for equality of the distribution functions.
As shown in Cambanis \textit{et al.}~\cite{Cambanis1981} $\vk{S}\stackrel{d}{=}
R_k\U$ is a spherically symmetric random vector with tractable
distributional properties. For
instance $(S_1 ,\ldots,S_m)^\top\stackrel{d}{=}R_m \U_m,m < k$ with positive
random radius $R_m$ such that
%
\begin{eqnarray}
\label{rm} R_m^2 \stackrel{d} {=}R^2_k
\kal{B}_{m/2,(k-m)/2},
\end{eqnarray}
with $\kal{B}_{m/2,(k-m)/2}$ independent of $R_k$. Equation \eqref{rm}
can be written iteratively as
%
\begin{eqnarray}
\label{it} R_{m}^2&\stackrel{d} {=}& R_{m+1}^2
\kal{B}_{m/2,1/2}, \qquad m=1 ,\ldots,k-1,
\end{eqnarray}
where $R_{m+1}^2 $ and $\kal{B}_{m/2,1/2}$ are independent. Note that
if $R_k^2$ is chi-square distributed with $k$ degrees of freedom (abbreviate
this by $R^2_k \sim\chi_k^2),$ then \eqref{it} holds for any $m\in
\N
$ with
$R_m^2\sim\chi^2_m$.

Another interesting result of Cambanis \textit{et al.}~\cite{Cambanis1981} is
that $ \boldsymbol{\mu}^\top\vk{S}\stackrel{d}{=}\sqrt {\boldsymbol{\mu
}^\top\boldsymbol{\mu}}
S_1$
for any $\boldsymbol{\mu}\in\R^k$. Consequently, the assumption that
$\Sigma
$ is a correlation matrix yields
\[
X_i\stackrel{d} {=}X_1\stackrel{d} {=}R_k
U_1, \qquad1 \le i \le k.
\]
We call a positive random variable $Z \sim F$ regularly varying at 0
with index $\gamma\in[0,\IF]$ if
%
\begin{eqnarray}
\label{eq:regv} \lim_{s \downarrow0} \frac{ F(st)}{F(s)}&=& t^\gamma \qquad
\forall t>0,
\end{eqnarray}
which is abbreviated as $Z \in RV_{\gamma}$ or $F \in RV_{\gamma}$.
Condition \eqref{eq:regv} is equivalent with $1/Z$ (or its survival function)
being regularly varying at infinity with index $-\gamma$. When $\gamma
=-\IF$, then the survival
function of $1/Z$ is called rapidly varying at infinity. See Jessen and
Mikosch~\cite{Jessen2006} or Omey and Segers~\cite{Omey2010} for
details on
regular variation.

Central for our results is an interesting fact discovered by Kabluchko
\cite{Kabluchko2011} pointing out the importance of the
incremental variance matrix (function) for the properties of the
Brown--Resnick process. Given a $k$-dimensional Gaussian random vector
$\X$ this $k\times
k$ matrix is denoted by $\Gamma=(\gamma_{ij})_{i,j\le k},$ where
$ \gamma_{ij} = \operatorname{\mathbf{Var}}\{X_i-X_j\}$. The
covariance matrix $\SIG$ of $\X$ is related to $\Gamma$ by
%
\begin{eqnarray}
\label{sig} \SIG=AA^\top= \bigl(\vks\vk{1}^\top+ \vk{1}
\vks^\top- \Gamma\bigr)/2, \qquad \vks=\bigl( \operatorname{
\mathbf{Var}}\{X_1\} ,\ldots ,\operatorname {\mathbf{Var}}
\{X_k\}\bigr)^\top.
\end{eqnarray}
If $\{Z(t),t\in T\}$ is a mean-zero Gaussian process with variance
function $\sigma^2(\cdot)$, we define similarly to the discrete case
the incremental variance function $\Gamma$ by
\[
\Gamma(t_1,t_2)= \operatorname{\mathbf{Var}}\bigl
\{Z(t_2)- Z(t_1)\bigr\},\qquad t_1,t_2
\in T.
\]
By Theorem~4.1 of Kabluchko~\cite{Kabluchko2011}, the stochastic process
%
\begin{eqnarray}
\label{kk} \eta_{\Gamma}(t) = \min_{i \ge1} \bigl\llvert
\Upsilon_i+ Z_i(t) \bigr\rrvert , \qquad t\in\R
\end{eqnarray}
is the limit of the minima of absolute values of independent Gaussian processes,
if additionally ${\Xi}_L = \sum_{i=1}^\IF\ve_{\Upsilon_i}$ is a
Poisson point process on $\R$ with points $\Upsilon_1,\Upsilon_2,\ldots$
and intensity measure given by the Lebesgue measure being further
independent of the Gaussian processes $\{Z_i(t), t\in\R\}, i\ge1$.
Here $\ve_x$ denotes the Dirac measure at $x$;
$\ve_x(B)=1$ if $x\in B\subset\R,$ and $\ve_x(B)=0$ when $x \notin B$.

In the sequel, for given $\boldsymbol{\theta}\in(0,\IF)^k, k\ge2$ and
$A,\SIG,\Gamma$ satisfying \eqref{sig} we write $\X\eqsim\mathfrak
{E}[\boldsymbol{\theta},\Gamma; H_k]$
if $\X\stackrel{d}{=}\ZW{ R_k A \vk{U}}, R_k\sim H_k$. We write simply
$\X\eqsim\mathfrak{{E}}[\Gamma; H_k]$ if the specification of
$\boldsymbol{\theta}$ is not
necessary for the stated result, meaning that
the result holds for any $\boldsymbol{\theta}\in(0,\IF)^k$.
Further, if
$R_k^2\sim\chi^2_k$ we write $\X\eqsim \mathit{Gauss}[ \Gamma]$, with $\X$
a mean-zero
Gaussian random vector with incremental variance matrix $\Gamma$.

\section{Minima of elliptical triangular arrays}\label{sec3}
Let $\X_{n}^{(i)}\stackrel{d}{=}R_k A_n \vk{U}, 1\le i \le n, n\ge
1$ be
$k$-dimensional independent elliptical random vectors, where the square matrix
$A_n$
is such that $\Sigma_n = A_nA_n^\top, n\ge1$ is a correlation
matrix. Next, we discuss the asymptotic behaviour of $\vk{L}_n =
(L_{n1} ,\ldots,
L_{nk})^\top, n\ge1$
defined by
\[
L_{nj} = \min_{1 \le i \le n} \bigl\llvert X_{nj}^{(i)}
\bigr\rrvert , \qquad j=1,\ldots,k, n\ge1.
\]
We have
\[
X_{nj}^{(i)} \stackrel{d} {=}X_{11}^{(1)}=:X_{11},
\qquad L_{nj} \stackrel{d} {=} L_{n1}, \qquad j=1,\ldots,k, 2
\le i\le n
\]
and $ \llvert  X_{11} \rrvert^2 \stackrel{d}{=}R^2_k \kal{B}_{1/2, (k-1)/2}.$

Next, we assume that $R_k \in RV_\gamma$ with index $\gamma\in
(0,1]$, which in view of {Lemma~\ref{lem:000}} implies $\llvert  X_{11}
\rrvert
\in RV_\gamma$;
note that the converse holds if $\gamma\in(0,1)$. Define a sequence
of constants $a_n,n\ge1$ by
%
\begin{eqnarray}
\label{an} \vk{P}\bigl\{a_n^{-1} \ge X_{11} > 0
\bigr\} &=& 1/n.
\end{eqnarray}
For such constants, we have the convergence in distribution ($n\to\IF$)
\begin{eqnarray*}
a_n L_{nj} &\stackrel{d} {\to}& \kal{L}_j
\sim\kal{G}_\gamma, \qquad j=1 ,\ldots,k,
\end{eqnarray*}
with distribution function $\kal{G}_\gamma$ given by
%
\begin{eqnarray}
\kal{G}_\gamma(x)&=& 1- \exp\bigl(-2x^\gamma\bigr), \qquad x> 0.
\end{eqnarray}
In view of Hashorva~\cite{Hashorva2011}, if $\SI_n$ has all
off-diagonal elements bounded by some constant $c\in(0,1)$, then
%
\begin{eqnarray}
\label{eq:allow} a_n \vk{L}_n \stackrel{d} {\to}\boldsymbol{\kal{L}} =(\kal{L}_1 ,\ldots ,\ZW{\kal {L}_k})^\top,
\qquad n\to\IF
\end{eqnarray}
holds with $\kal{L}_1 ,\ldots,\ZW{\kal{L}_k}$ being mutually independent.
By allowing the off-diagonal elements of $\SI_n$ to converge to 1 as
$n\to\IF$ with a certain speed, it is possible that
the random vector $\boldsymbol{\kal{L}}$ has dependent components. If
$H_i,i\le k$ is the distribution function of $\ZW{R_i}$ in \eqref{it}
it turns out that $\kal{R}_m, m\le k-1$ with distribution function
%
\begin{eqnarray}
\label{barh} \kal{H}_m(z)& = & \int_0^z
\frac{1}{r \vk{E}\{1/R_{m+1}\} } \,\mathrm{d}H_{m+1}(r), \qquad z>0
\end{eqnarray}
determine the distribution function of $\boldsymbol{\kal{L}}$ (assuming $\vk
{E}\{1/R_k\}< \IF$). For the derivation of this result, we shall define
an elliptical random vector $\vk{Z}^{K;j}\stackrel{d}{=}\kal
{R}_{m-1}\Gamma_{m,K}\vk{U}_m$ with
\begin{eqnarray*}
&&\Gamma_{m,K}(\Gamma_{m,K})^\top=\bigl( \vk{1}
\Gamma_{K_j,J}^\top+\Gamma_{K_j,J} \vk{1}^\top-
\Gamma_{K_j,K_j}\bigr)/2, \\
&& \quad \vk{1}=(1 ,\ldots,1)^\top\in
\R^{m-1}, K_j=K\setminus J, J=\{j\},
\end{eqnarray*}
where $K\subset\{ 1 ,\ldots,k\}$ has $m\ge2$ elements, and $\Gamma$ is
the matrix in \eqref{eq:reis}.

\begin{theo}\label{th:kabl:0} Let $\X_{n}^{(i)}, 1\le i \le n, n\ge
1$ be a
triangular array of $k$-dimensional elliptical random vectors with
correlation matrices $\Sigma_n,n\ge1$ as above, and $R_k\sim H_k$.
Suppose that $\llvert  X_{11}^{(1)} \rrvert \in RV_\gamma,\gamma\in
(0,1]$ and $\vk{E}\{1/R_k\}< \IF$.

If condition \eqref{eq:reis} is satisfied for $c_n = 2 a_n^2$ with
$a_n$ determined by \eqref{an}, then \eqref{eq:allow} holds and
for all $\vk{x}\in(0,\IF)^k$
%
\begin{eqnarray}
\label{int} &&\vk{P}\{\boldsymbol{\kal{L}}> \vk{x}\} = \exp \Biggl( \ZW{\sum
_{m=1}^{k}(-1)^{m}}\sum
_{\llvert  K \rrvert =m} \int_{-x_j^{\gamma
}}^{x_j^{\gamma}} \vk{P}\bigl
\{\bigl| \sign(y) \llvert y \rrvert^{1/\gamma} +Z^{K;j}_i \bigr|
\nonumber
\\[-8pt]
\\[-8pt]
&&\hphantom{\vk{P}\{\boldsymbol{\kal{L}}> \vk{x}\} = \exp \Biggl( \ZW{\sum
_{m=1}^{k}(-1)^{m}}\sum
_{\llvert  K \rrvert =m} \int_{-x_j^{\gamma
}}^{x_j^{\gamma}} \vk{P}
\{}\le x_i, i\in K\setminus\{j\}, j\in K\bigr\}\,\mathrm{d} y \Biggr),
\nonumber
\end{eqnarray}
where the summation above runs over all non-empty index sets $K$ with
$\llvert  K \rrvert =m$ elements and $j$ is some index in $K$. Set the
integral in
\eqref{int} equal to $2 x_j^\gamma$ if $K=\{j\}$.
\end{theo}

\begin{remarks*}
\begin{enumerate}[(b)]
\item[(a)] The result of {Theorem~\ref{th:kabl:0}} can be extended
for $\Gamma$ with off-diagonal elements equal to 0. For instance
when $\Gamma= \vk{0} \vk{0}^\top$ with $\vk{0} = (0 ,\ldots
,0)^\top
$, then it follows that
\begin{eqnarray*}
\vk{P}\{\boldsymbol{\kal{L}}> \vk{x}\} &=& 1- \kal{G}_{\gamma}\Bigl(
\min_{1 \le i\le k} x_i\Bigr), \qquad\vk{x}\in (0,
\IF)^{k}.
\end{eqnarray*}
\item[(b)] In view of \eqref{int} the random vector $(\kal{L}_d,\kal{L}_l),
d\not=l$ has joint distribution function depending on the element
$\gamma_{dl}$ of $\Gamma$.
\end{enumerate}
\end{remarks*}
\begin{example}\label{exa1} Let $\X_{n}^{(i)}, 1 \le i\le n, n\ge1$ be a
triangular array of $k$-dimensional mean-zero Gaussian random vectors with
covariance matrix $\SI_n,n\ge1$. Since $\kal{R}^2_m \sim\chi^2_m,m\le k$, then
$a_n$ defined by \eqref{an} satisfies
\[
a_n= \bigl(1+\mathrm{o}(1)\bigr)\frac{n}{\sqrt{2 \uppi}}, \qquad n\to\IF.
\]
Hence, when \eqref{eq:reis} is valid with $c_n = 2 a_n^2$, then \eqref
{int} holds with $\vk{Z}^{K;j}$ a mean-zero Gaussian random vector with
covariance matrix $\Gamma_{m,K}(\Gamma_{m,K})^\top$.
\end{example}

Next, we extend {Theorem~\ref{th:kabl:0}} imposing a smoothness
assumption on
$R_k$, namely that \eqref{it} holds also for $m=k$.

\begin{theo}\label{T2}
Under the assumptions and notation of {Theorem~\ref{th:kabl:0}}, if further
\eqref{it} holds for $m=k$ with $R_{k+1}\sim H_{k+1}$, then 
%
\begin{eqnarray}
\label{eq:genralise:p} \vk{P}\{\boldsymbol{\kal{L}}> \vk{x}\} &=& \exp \biggl(- \int
_{\R} \vk{P}\bigl\{\exists i \le k\dvt  \bigl\llvert \sign (y)
\llvert y \rrvert^{1/\gamma}+Z_{i} \bigr\rrvert \le
x_i \bigr\} \,\mathrm{d} y \biggr), \qquad\vk{x}\in(0,\IF)^{k}, \qquad\quad
\end{eqnarray}
with $\vk{Z} \eqsim\mathfrak{E}[\Gamma; \kal{H}_{k}]$ and $\kal
{H}_{k}$ defined by \eqref{barh}.
\end{theo}

\begin{remark*}
The assumption \eqref{it} is satisfied for $m=k$, provided that $\X_n^{(i)},i\le n$ is a subvector of an elliptical random vector, see
Cambanis \textit{et al.}~\cite{Cambanis1981}. In particular, it holds if $R_k
\stackrel{d}{=}S \tilde{R}_k$ with $S$ a positive random variable
independent of $\tilde{R}_k^2 \sim\chi^2_{k}$.
\end{remark*}
\begin{example} Let $\X_{n}^{(i)}, 1\le i\le n, n\ge1$ be as in
Example~\ref{exa1}. Next, define
\[
\Y_n^{(i)} = S_n\X_n^{(i)},
\qquad1 \le i\le n, n\ge1,
\]
%
with
$S,S_n,n\ge1$ independent positive random variables with distribution
function $F$ being further independent of $\X_{n}^{(i)}, 1 \le i\le n$.
If $F \in RV_\gamma, \gamma\in(0,1]$, then 
by {Lemma~\ref{lem:000}} $\llvert  Y_{n1}^{(1)} \rrvert \in RV_{\gamma}$.
Define constants
$a_n,n\ge1$ such that $\vk{P}\{0< S X_{11}(1)\le1/a_n\} =1/n$ holds for
all large $n$. If further \eqref{eq:reis} is satisfied with $c_n=2a_n^2$,
then \eqref{eq:genralise:p} holds. Note in passing that $\kal{H}_{k}$
satisfies \eqref{barh} with
$R^2_{k+1} \sim\chi^2_{k+1},R_{k+1}>0$.
\end{example}
\section{Maxima of elliptical triangular arrays}\label{sec4}
With the same notation as above we consider again the triangular array
$\X_{n}^{(i)}, 1\le i \le n, n\ge1$ of $k$-dimensional independent
elliptical random vectors with stochastic
representation \eqref{eq:el} and $\Sigma_n = A_nA_n^\top, n\ge1$
given correlation matrices. Define the componentwise maxima $ \vk{M}_n
= (M_{n1}
,\ldots,M_{nk})^\top$ by
\[
M_{nj} = \max_{1 \le i \le n} X_{nj}^{(i)}, \qquad
j=1 ,\ldots,k, n\ge1.
\]
The asymptotic behaviour of the maxima of elliptical triangular arrays
is discussed in Hashorva~\cite{Hashorva2006} assuming that
the random radius $R_k$ has distribution function $H_k$ in the Gumbel
MDA. A canonical example of such triangular arrays is that of the
Gaussian arrays for which
the limit distribution of the maxima is the H\"usler--Reiss copula which
is a particular case of the Brown--Resnick copula.
When $H_k$ is in the Weibull MDA the limit distribution of the maxima
is a max-infinitely divisible distribution function, see Hashorva \cite
{Hashorva2005}.

We reconsider the findings of the aforementioned papers showing novel
representations of the limit distributions given in terms of the
distribution of the maxima of some point processes shifted by
elliptical random
vectors.
For the derivation of the next results, we impose asymptotic
assumptions on either the marginal distribution functions or on the
associated random radius
$R_k$, which is of some interest for statistical applications where
some data might be missing, or some component of the random vector
might be
unobservable, and therefore the random radius itself cannot be estimated.

\subsection{Gumbel max-domain of attraction}\label{sec4.1}
The main assumption in this section is that the marginal distribution
functions of the elliptical triangular array are in the Gumbel MDA.
A univariate distribution function $G$ is in the Gumbel MDA
(abbreviated $G\in\mathrm{GMDA} (w)$) if for any $x\in\R$
%
\begin{eqnarray}
\label{eq:gumb1} \lim_{t\uparrow\omega} \frac{1-
G(t+ x/w(t))}{1-
G(t)}=\exp(-x), \qquad\omega=
\sup\bigl\{t\dvt  G(t)<1\bigr\},
\end{eqnarray}
with $w(\cdot)$ some positive scaling function. If $\omega=\IF$, an
important property for the distribution function $G$ satisfying \eqref
{eq:gumb1} is a key finding of Davis and Resnick~\cite{Davis1988},
namely by Proposition~1.1 therein (see also Embrechts \textit{et al.} \cite
{Embrechts1997}, page~586) for any $\mu\in\R, \tau>1$ we have
%
\begin{eqnarray}
\label{Davis} \lim_{x \to\infty} \bigl(x w(x)\bigr)^\mu
\frac{1- G(\tau x)}{1- G(x)}&=& 0.
\end{eqnarray}
Indeed \eqref{Davis}, which we refer to as the \textit{Davis--Resnick
tail property} is crucial for several asymptotic approximations.

\begin{theo}\label{theo:main} Let $R\sim H_k, \X_n^{(i)}, 1\le i\le
n, \SI_n,n\ge1$ be as in {Theorem~\ref{th:kabl:0}}.
If either $G\in\mathrm{GMDA}(w)$ or $H_k\in\mathrm{GMDA}(w)$ and
condition \eqref{eq:reis} is satisfied with
%
\begin{eqnarray}
\label{eq:main:cond}\label{ab} c_n = 2\frac{b_n}{a_n}, \qquad
b_n = G^{-1}(1- 1/n),  a_n =
1/w(b_n),  n>1,
\end{eqnarray}
then for any $\vk{x}\in\R^k$ and $\vk{Z}\eqsim \mathit{Gauss}[\Gamma]$ we have
%
\begin{eqnarray}
&&\lim_{n \to\infty} \vk{P}\bigl\{(\vk{M}_n - b_n \vk{1})
/a_n \le\vk {x}\bigr\}  \nonumber
\\[-8pt]
\\[-8pt]&& \quad = Q_{\Gamma}(\vk{x})= \exp \biggl(-
\int_{\R} \vk{P}\{\exists i \le k\dvt  Z_i >
x_i-y + \theta_i/2 \} \exp(-y)\,\mathrm{d}y \biggr),
\nonumber
\end{eqnarray}
where $\theta_i=\operatorname{\mathbf{Var}}\{Z_i\},i \le k$.
\end{theo}
Since the above result holds for Gaussian triangular arrays with
scaling function \mbox{$w(x)=x$}, the distribution function $Q_{\Gamma}$ is
the multivariate max-stable H\"usler--Reiss distribution function. For a
particular choice of a Gaussian process $\{Z(t),t\in\R\}$
this distribution has the Brown--Resncik copula; in fact it can be
directly defined by Brown--Resnick processes $\beta_{\kal{R};\Gamma}$
with independent
Gaussian points $\xi_i(t):=Z_i(t)- \sigma^2(t)/2, i\ge1$ given as
%
\begin{eqnarray}
\label{eq:BRBR} \beta_{\kal{R};\Gamma}(t) = \max_{i \ge1} \bigl[
\Upsilon_i+ \xi_i(t)\bigr], \qquad t\in\R.
\end{eqnarray}
Here $ \Xi=\sum_{i=1}^\IF\ve_{\Upsilon_i}$ is a Poisson point
process with intensity measure $\exp(-x)\,\mathrm{d}x$ being independent of
$\{Z_i(t), t\in\R\}, i\ge1$.
In view of our result, the Brown--Resnick process with Gaussian points
does not depend on the
variance function, which is already established in Theorem~2.1 of
Kabluchko \textit{et al.}~\cite{Kabluchko2009}.

\COM{
The explicit form of the limit distribution function $Q_{\Gamma}$
suggests a novel class of max-stable multivariate distributions, i.e.,
for some
random vector $\vk{Z}= B \kal{R}_k\boldsymbol{\kal{U}}$ where $\kal{R}_k$
is a positive random variable independent of $\boldsymbol{\kal{U}}$ with
$R_k^2 \sim\chi_k^2$ we may define
%
\begin{eqnarray}
\label{Tnew} Q_{B,\boldsymbol{\kal{U}}}(\vk{x})&=& \exp \biggl(- \int_{\R}
\vk{P}\{ \exists i \le k: Z_i> x_i-y +
\theta_i/2 \} \exp(-y)\,\mathrm{d}y \biggr) \qquad \forall\vk{x}\in
\R^k,
\end{eqnarray}
where $B$ be is a $k\times k $ matrix given by
%
\begin{eqnarray}
\label{m:B} BB^\top&=& \bigl(\vk{1} \boldsymbol{\theta}^\top+
\boldsymbol{\theta }\vk{1}^\top- \Gamma\bigr)/2, \qquad\boldsymbol{
\theta} \in(0,\IF)^k.
\end{eqnarray}
%
The random vector $\boldsymbol{\kal{U}}$ needs to obey certain conditions
having components $\kal{U}_i$ with upper endpoints say equal to
$\omega\in(0,\IF]$ and $\vk{E}\{\exp( R_k U_i)\}< \IF, i\le n$.}

\subsection{Weibull max-domain of attraction}\label{sec4.2}


The unit Weibull distribution with index $\alpha\in(0,\IF)$ is $\Psi_\alpha(x)=\exp(- \llvert  x \rrvert^\alpha), x<0$.
In view of Hashorva and Pakes~\cite{Hashorva2010} the distribution
function $G$ is in the Weibull MDA if $H_k$ is in the Weibull MDA. We
assume for simplicity that $H_k$ has upper endpoint
equal to 1. By definition, $H_k$ is in the MDA of $\Psi_\alpha$ (for
short $H_k\in\mathrm{WMDA} (\alpha)$) if for any $x\in(0,\IF)$
%
\begin{eqnarray}
\label{eq:kela} \lim_{n \to\infty} H^n_k\bigl(1-a(n)x
\bigr) &=& \Psi_{\alpha}(x) , \qquad a_n = 1-H^{-1}_k(1-
1/n).
\end{eqnarray}

If $H_k\in\mathrm{WMDA} (\alpha)$, with some index $\alpha\in(0,\IF
)$ and $H_k$ has upper endpoint equal to~1, then by Theorem~2.1 in
Hashorva~\cite{Hashorva2008}
%
\begin{eqnarray}
\label{happy} \lim_{n \to\infty} \vk{P}\bigl\{(\vk{M}_n - \vk{1})
/a_n \le\vk{x}\bigr\} &=& \widetilde {\kal{Q}_{\Gamma, \alpha}}(
\vk{x}) \qquad\forall\vk{x}\in(-\IF,0)^k,
\end{eqnarray}
with $\widetilde{\kal{Q}_{\Gamma, \alpha}}$ a max-infinitely
divisible distribution function, provided that \eqref{eq:reis} holds
with $c_n = 2/a_n, a_n = 1-
G^{-1}(1- 1/n), n>1$.

In the next theorem, we show that \eqref{happy} holds if either $G$ or
$H_k$ is in the Weibull MDA. Furthermore,
we give a new representation for the limit distribution function
$\widetilde{\kal{Q}_{\Gamma, \alpha}}$.

\begin{theo}\label{theo:Weib} Let $R\sim H_k, \X_n^{(i)}, 1 \le i\le
n, \SI_n,n\ge1$ be as in {Theorem~\ref{th:kabl:0}}, and assume that $G$
has upper endpoint~1. If either $G \in\mathrm{WMDA} (\alpha+ (k-1)/2),
$ or $H_k\in\mathrm{WMDA} (\alpha)$, with $\alpha\in(0,\IF)$, then
\eqref{happy} holds where
%
\begin{eqnarray}
\label{reff} 
\widetilde{\kal{Q}_{\Gamma, \alpha}}(
\vk{x})&=& \exp \biggl(- \int_{0}^\IF\vk{P}\{
\exists i \le k\dvt  \sqrt{2 y} Z_i > x_i+y +
\theta_i/2 \} \,\mathrm{d} y^{\alpha+ (k-1)/2} \biggr),
\end{eqnarray}
with $\vk{Z}\eqsim\mathfrak{{E}}[\Gamma; H_k], \vks\in(0,\IF)^k$
and $\widetilde{\kal
{H}}_\alpha$ the distribution function of $\widetilde{ \kal
{R}_{\alpha}}>0$ which satisfies
$\widetilde{\kal{R}_\alpha}^2 \stackrel{d}{=}\kal{B}_{k/2, \alpha}$.
\end{theo}
We remark that $\widetilde{\kal{Q}_{\Gamma, \alpha}}$ has Weibull
marginal distributions $\Psi_{\alpha+ (k-1)/2}$.
It follows from our result that $\widetilde{\kal{Q}_{\Gamma, \alpha
}}$ is determined by $\Gamma$ and $\alpha$ but not by the vector
$\boldsymbol
{\theta}$,
and further $\widetilde{\kal{Q}_{\Gamma, \alpha}}$ is not a
max-stable distribution function; clearly, it is a max-infinitely
divisible distribution function.

\section{Results for spherical processes}\label{sec5}
It is well-known that spherical random sequences are mixtures of
Gaussian random sequences.
Specifically, if the random variables $X_i,i\ge1$ with some common
non-degenerate distribution function $G$ are such that $(X_1 ,\ldots,
X_k)$ is centered and spherically distributed
for any $k\ge1$, then $X_i\stackrel{d}{=}S X_i^*, i\ge1$ with
$X_i^*,i\ge
1$ is a sequence of independent standard
Gaussian random variables being further independent of $S>0$.
Consequently, a spherical random process $\{ X(t), t\in\R\}$ such
that $X(t)$ has distribution function $G$
for any $t\in\R$
can be expressed as $\{X(t)= SY(t),t\in\R\}$ with $Y(t)$ a mean-zero
Gaussian process and $S$ a positive random variable independent of
$\{Y(t),t\in\R
\}$; see Theorem~7.4.4 in Bogachev~\cite{Bogachev1998} for a general
result on spherically symmetric measures. We note in passing that $\{
X(t), t\in T\}$
is a particular instance of Gaussian processes with random variance,
see H\"usler \textit{et al.}~\cite{H2011b} for recent results on extremes of
those processes.

We shall discuss first the asymptotic behaviour of the maximum of
independent spherical processes.
Then we shall briefly investigate the asymptotics of the minima of
absolute values of those processes.

\textit{Model A}: Assume that $S$ has an infinite upper endpoint such that
for given constants $\alpha_1\in\R$ and $C_1,L_1, p_1 \in(0,\IF)$
%
\begin{eqnarray}
\label{Weibullian} \vk{P}\{S> x\} &=&\bigl(1+\mathrm{o}(1)\bigr) C_1
x^{\alpha_1} \exp\bigl(-L_1 x^{p_1}\bigr), \qquad x\to
\IF
\end{eqnarray}
is valid. We abbreviate \eqref{Weibullian} as $S\in\kal
{W}(C_1,\alpha_1, L_1,p_1)$.

\textit{Model B}:
Consider $S$ with upper endpoint equal to 1 such that
%
\begin{eqnarray}
\label{eq:SA} \lim_{u \to\infty} \frac{\vk{P}\{S > 1- x/u\} }{\vk{P}\{S> 1-
1/u\} }&=& x^\gamma, \qquad
x\in(0,\IF),
\end{eqnarray}
with $\gamma\in[0,\IF)$ some constant.

Since for $S=1$ almost surely, the spherical process is simply a
Gaussian one
(which is covered by Model B for $\gamma=0$) intuitively, we expect
that under the Model B the maximum of independent elliptical processes will
behave asymptotically as the maximum of independent Gaussian processes.
This intuition is confirmed by {Theorem~\ref{th:last}} below.
In fact, it turns out that the limit process of the maximum of
independent spherical processes is in both models the
Brown--Resnick process.
Next, if $\Gamma(\cdot,\cdot)$ is a negative definite kernel in $\R^2$ we define as previously the Brown--Resnick stochastic process with
Gaussian points
as
%
\begin{eqnarray}
\label{def:zeta} \beta_{\kal{R};\Gamma}(t) = \max_{i\ge1} \bigl(
\Upsilon_i+ Z_i(t) -\sigma^2(t)/2 \bigr),
\qquad t\in T \subset\R,
\end{eqnarray}
with $\{Z_i(t), t\in T\}$ independent Gaussian processes with
incremental variance function~$\Gamma$, variance function $\sigma^2(\cdot)$ being further
independent of the point process $ \Xi$ with points $\Upsilon_i,i\ge
1$ appearing in \eqref{eq:BRBR}. For simplicity, we deal below with
the case
$T=\R$ establishing weak convergence of finite-dimensional
distributions (denoted below as $ \Longrightarrow$).

\begin{theo}\label{th:last} Let $\{Y_{ni}(t), t\in\R\}, 1\le i\le n,
n\ge1$
be independent Gaussian processes with mean-zero, unit variance
function and
correlation function $\rho_n(s,t), s,t\in\R$. Let $S,S_{ni},n\ge1$ be
independent and identically distributed positive random variables. Set
$ \{X_{ni}(t) = S_{ni} Y_{ni}(t),t\in\R\}, n\ge1$, and let
$G$ be the distribution function of $X_{11}(1)$. Suppose that
%
\begin{eqnarray}
\label{rho:fish} \lim_{u \to\infty} c_n 
\bigl(1-
\rho_n(t_1,t_2) \bigr) &= &
\Gamma(t_1,t_2)\in(0,\IF), \qquad t_1
\not=t_2 \in\R,
\end{eqnarray}
where $c_n=2 b_n/a_n$ and $a_n=1/w(b_n), b_n= G^{-1}(1- 1/n)$ with
$G^{-1}$ the inverse of $G$. 
\begin{enumerate}[(B)]
\item[(A)] If \eqref{Weibullian} holds, then as $n\to\IF$
%
\begin{eqnarray}
\label{claim:M} \frac{1}{a_n} \Bigl[\max_{1 \le i \le n}X_{ni}(t) -
b_n\Bigr]&\Longrightarrow& \beta_{\kal{R};\Gamma}(t), \qquad t\in\R,
\end{eqnarray}
where $\Longrightarrow$ means the weak convergence of the
finite-dimensional distributions, and
\begin{eqnarray*}
&&\frac{b_n}{a_n} =\bigl(1+\mathrm{o}(1)\bigr) \frac{2p_1\ln n}{2+p_1} , \qquad
b_n= \bigl(1+\mathrm{o}(1)\bigr) \biggl( \frac{\ln n}{L_{1}A^{-p_{1}}+A^{2}/2}
\biggr)^{(2+p_1)/(2 p_1)}, \\
&& \quad  A=(p_1L_1)^{1/(2+p_1)}.
\end{eqnarray*}
\item[(B)] If \eqref{eq:SA} holds with $\gamma\in[0,\IF)$, then \eqref
{claim:M} is satisfied and $\lim_{n \to\infty} b_n/\sqrt{2 \ln
n}=\lim_{n \to\infty}
a_n\sqrt{2 \ln
n}=1$.
\end{enumerate}
\end{theo}

Next, we discuss the asymptotic behaviour of the minimum of absolute
values in the framework of independent spherical processes.

\begin{theo}\label{th:last:2} Let $\{Y_{ni}(t), Z_i(t), t\in\R\},
1\le i\le
n, n\ge1 $ be as in {Theorem~\ref{th:last}}, and let $\{S_{ni}(t),t\in
\R
\},
n\ge1$ be
independent copies of $\{S(t),t\in\R\}$, being further independent of
the Gaussian processes. Define the spherical processes $ \{X_{ni}(t) = S_{ni}(t)
Y_{ni}(t),t\in\R\}, n\ge1$, and suppose that $S(t)>\kappa, t\in\R$
almost surely for some positive constant $\kappa$. If $a_n= n/\sqrt{2
\uppi}$ and \eqref{rho:fish} holds with $c_n=2 a_n^2$,
then as $n\to\IF$
%
\begin{eqnarray}
\label{howpray:toJ} \min_{1 \le i \le n} a_n \bigl\llvert
X_{ni} (t) \bigr\rrvert &\Longrightarrow& \min_{i\ge1}
S_i(t)\bigl| \Upsilon_i+ Z_i(t) \bigr|=
\zeta_{\Gamma
,S}(t), \qquad t\in\R,
\end{eqnarray}
where $\Upsilon_i,i\ge1$ are the points of $\Xi$ defined in (\ref{eq:BRBR})
being independent of both $Z_i(t),S_i(t),\break t\in\R, i\ge1$.
\end{theo}

\begin{remarks*}
\begin{enumerate}[(b)]
\item[(a)] In {Theorem~\ref{th:last:2}}, we can relax the assumption
that $S(t)$ is bounded from below by assuming instead $\vk{E}\{
[S(t)]^{-1-\ve}\}<
\IF$ for some $\ve>0$.

\item[(b)] The process $\{\zeta_{\Gamma, S}(t), t\in\R\}$ is defined by
$\Gamma$ and $\{S(t),t\in\R\}$ but does not depend on the variance function
$\sigma^2(\cdot)$. The processes $\zeta_{\Gamma,1}$ appears first
in Penrose~\cite{Penrose1991} and recently in Kabluchko \cite
{Kabluchko2011}. We refer to
$\{\eta_{\Gamma,S}(t),t\in\R\}$ as Penrose--Kabluchko process.
\end{enumerate}
\end{remarks*}
\section{Further results and proofs}\label{sec6}

\begin{lem}\label{lem:000} Let $\X\stackrel{d}{=}R A \vk{U}$ be an
elliptical
random vector in $\R^k,k\ge2$ with $A$ such that $AA^\top$ is a
positive definite
correlation matrix and $R>0$.
\begin{enumerate}[(b)]
\item[(a)] If for some $\gamma\in[0,\IF]$ we have $R\in RV_{\gamma}$, then
$\llvert  X_1 \rrvert \in RV_{\gamma^*}$ with $\gamma^* = \min( \gamma,1)$.

Conversely, if $\llvert  X_1 \rrvert \in RV_{\gamma^*}$ with $\gamma^*\in(0,1)$, then $R\in RV_{\gamma^*}$.

\item[(b)] If $\vk{E}\{R^{-1- \ve}\} < \IF$ for some $\ve>0$, then $\llvert  X_1
\rrvert \in RV_1$.\vadjust{\goodbreak}
\end{enumerate}
\end{lem}

\begin{pf} (a) If $\gamma\in[0,\IF)$ the
proof follows from
Theorem~4.1 in Hashorva~\cite{Hashorva2011}.
When $\gamma=\IF$, then $1/R$ is rapidly varying at infinity. Hence
from Theorem~5.4.1 of de Haan and Ferreira~\cite{De2006}
$\vk{E}\{R^{-p}\}< \IF$ for any $p\in(0,\IF)$, and thus the claim follows
once the statement (b) is proved.
Statement (b) can be directly established by applying Breiman's lemma
(see Breiman~\cite{Breiman1965}, Davis and Mikosch~\cite{Davis2008}),
and thus the proof is complete.
\end{pf}

\begin{pf*}{Proof of Theorem~\ref{th:kabl:0}}
By the relation between the
minima and maxima,
in view of Lemma~4.1.3 in Falk \textit{et al.}~\cite{Falk2010}
the proof follows if
%
\begin{eqnarray}
\label{show} \lim_{n \to\infty} n \vk{P}\bigl\{a_n \llvert
X_{ni} \rrvert \le x_i, i\in K\bigr\} &=&
L_K(\vk{x}_K), \qquad\vk{x}\in(0,\IF)^k
\end{eqnarray}
holds for any non-empty index set $K\subset\{1 ,\ldots,k\}$ with
$m\ge
2$ elements, and $L_K(\cdot)$ some right-continuous functions.
In the sequel, we write simply $\X_n$ instead of $\X_n^{(1)}$; the
subvector $(\X_n)_K$ is an elliptical random vector with associated
random radius
$R_m\sim H_m$ satisfying \eqref{it}. By {Lemma~\ref{lem:000}}, $H_k
\in
RV_\gamma, \gamma\in(0,1]$ implies $H_m\in RV_\gamma, 1 \le m\le k-1$.
Consequently, it suffices to show \eqref{show} for the case $m=k$.
Since the distribution function of $\X_n$ depends on $\Sigma_n$ and
not on $A_n$, and further
$\Sigma_n$ is positive definite, we can assume that $A_n$ is a lower
triangular matrix. Define $q_n(y) = y/a_n,y \in\R$ and recall that $G$
denotes the distribution function of $X_{11}$.
It follows that conditioning on $X_{nk}=q_n(y)$ with $y\not=0$ such
that $G(\llvert  y \rrvert /a_n)\in(0,1),n \ge1$
we have the stochastic representation (set $I=\{1 ,\ldots,k-1\}, J=\{
k\}$)
%
\begin{eqnarray}
\label{alves:0} (\X_{n})_I | X_{nk}=q_n(y)
&\stackrel{d} {=}& R_{y,n,k-1} B_{nk} \vk{U}_{k-1} + (
\Sigma_{n})_{IJ} q_n(y), \qquad n\ge1,
\end{eqnarray}
where $B_{nk}$ is a lower triangular matrix satisfying $B_{nk}
B_{nk}^\top=\ZW{(\Sigma_n)_{II}}- (\Sigma_{n})_{IJ}(\Sigma_{n})_{JI} .$
In view of Cambanis \textit{et al.}~\cite{Cambanis1981}, $\vk{U}_{k-1}$ is
independent of $R_{y,n,k-1}, n\ge1$ which has survival function
$\overline{Q}_{y,n,k-1}$ given
by
%
\begin{eqnarray}
\label{eq:surv} &&\overline{Q}_{y,n,k-1}(z)=\frac{\int_{((y/a_n)^2+
z^2)^{1/2}}^{\omega} (r^2- (y/a_n)^2)^{(k-1)/2- 1}\ZW{ r^{-k+2}}\,\mathrm{d}
H_k(r)} {
\int_{y/a_n}^{\omega} (r^2- (y/a_n)^2)^{(k-1)/2- 1}\ZW{r^{-k+2}} \,\mathrm{d} H_k(r)},
\nonumber
\\[-8pt]
\\[-8pt]
&& \quad  z\in
 \bigl(0, \sqrt{\omega^2- y^2/a_n^2}
 \bigr).
\nonumber
\end{eqnarray}
Clearly, $\lim_{n \to\infty} a_n=\IF$ and the monotone convergence theorem
implies the convergence in distribution
\[
R_{y,n,k-1} \stackrel{d} {\to}\kal{R}_{k-1}, \qquad n\to\IF,
\]
where $\kal{R}_{k-1} \sim\kal{H}_{k-1}$ with
%
\begin{eqnarray}
\label{qN} \kal{H}_{k-1}(z)
=1-
\frac{\int_{z}^{\omega} r^{-1}\,\mathrm{d}H_k(r)} { \vk{E}\{1/R_k\}}, \qquad z\in(0, \omega).
\end{eqnarray}
In view of relation \eqref{rm} and since for any integer $m\ge2$ we
have $\vk{E}\{1/\kal{B}_{m/2, (k-m)/2}\}< \IF$ the assumption $\vk
{E}\{1/R_k\}<
\IF$
implies
$\vk{E}\{1/R_{m}\}< \IF$. Hence, the above convergence holds also for the
omitted case $k=m$.\vadjust{\goodbreak}
Next, by \eqref{eq:reis} and the fact that $B_{nk}B_{nk}^\top$ (and
not the matrix $B_{nk}$) defines the conditional distribution in
\eqref{alves:0}
we can choose $B_{nk}$ such that $\lim_{n \to\infty} a_n B_{nk}=
B_k$ with
\[
B_kB_k^\top= \bigl(\vk{1} \boldsymbol{
\theta}^\top+ \boldsymbol{\theta }\vk{1}^\top-
\Gamma_{II}\bigr)/2, \qquad\boldsymbol{\theta}= \Gamma_{IJ} .
\]
Hence, for any $\vk{x}\in(0,\IF)^k$ utilising further \eqref
{alves:0} and the fact that $G$ is symmetric about 0 we obtain (set
$G_n( y)=G(y/a_n),n\ge1$ %
and $K=\{1 ,\ldots,k\}$)
\begin{eqnarray*}
&&
\vk{P}\bigl\{a_n \llvert X_{ni}
\rrvert \le x_i, \forall i=1 ,\ldots,k\bigr\}
\\[-2pt]
&& \quad = \int_{\R} \vk{P}\bigl\{a_n \llvert
X_{ni} \rrvert \le x_i, \forall i\in I|
X_{nk}=y\bigr\} \,\mathrm{d} G(y)
\\[-2pt]
&& \quad = \int_{-x_k}^{x_k} \vk{P}\bigl\{a_n
\llvert X_{ni} \rrvert \le x_i, \forall i\in I|
X_{nk}=y/a_n\bigr\} \,\mathrm{d} G_n(y)
\\[-2pt]
&& \quad = \int_{0}^{x_k} \bigl[\vk{P}\bigl
\{a_n \llvert X_{ni} \rrvert \le x_i, \forall
i\in I| X_{nk}=y/a_n\bigr\} 
\\[-2pt]
&& \hphantom{= \int_{0}^{x_k} \bigl[}\quad {}+\vk{P}\bigl
\{a_n \llvert X_{ni} \rrvert \le x_i, \forall
i\in I| X_{nk}=-y/a_n\bigr\} \bigr]\,\mathrm{d} G_n(y)
\\[-2pt]
&& \quad =\int_{0}^{x_k} \bigl[\vk{P}\bigl
\{a_n \llvert Z_{ni}+ d_{ni} y/a_n
\rrvert \le x_i, i\in I\bigr\} +\vk{P}\bigl\{a_n \llvert
Z_{ni}- d_{ni} y/a_n \rrvert \le
x_i,\forall i\in I \bigr\} \bigr]\,\mathrm{d} G_n(y),
\end{eqnarray*}
with $\vk{Z}_{n}= R_{y,n,k-1} B_{nk}\vk{U}_{k-1}$ and $d_{ni}$ the
$i$th component of $(\Sigma_{n})_{IJ}$. By the construction we have
the convergence in distribution $(n\to\IF)$
\[
R_{y,n,k-1} (a_nB_{nk}) \vk{U}_{k-1}
\stackrel{d} {\to}\kal{R}_{k-1} B_k \vk {U}_{k-1}=:(Z_1
,\ldots,Z_{k-1})^\top.
\]
Further, by the regular variation at 0 of the distribution function of
$\llvert  X_{11} \rrvert $, the fact that $X_{11}$ is symmetric about 0,
and the choice of $a_n,n\ge1$ we have
%
\begin{eqnarray}
\label{quantitative} \lim_{n \to\infty} n \bigl[G_n( t)-
G_n( s)\bigr]&=& t^\gamma- s^\gamma \qquad\forall
s,t\in(0, \IF).
\end{eqnarray}
Consequently, since $\lim_{n \to\infty} d_{ni}= 1$
\begin{eqnarray*}
&&\lim_{n \to\infty} n \vk{P}\bigl\{a_n \llvert X_{ni}
\rrvert \le x_i, \forall i=1 ,\ldots,k\bigr\} \\[-2pt]
&& \quad = \int
_{0}^{x_k} \bigl[\vk{P}\bigl\{\llvert
Z_i+ y \rrvert \le x_i, i\in I \bigr\} \,\mathrm{d}
y^\gamma+ \vk{P}\bigl\{\llvert Z_i- y \rrvert \le
x_i, i\in I \bigr\} \bigr]\,\mathrm{d} y^\gamma
\\[-2pt]
&& \quad = \int_{0}^{x_k^{\gamma}} \bigl[ \vk{P}\bigl\{\bigl
\llvert Z_i+ y^{1/\gamma
} \bigr\rrvert \le x_i, i\in
I \bigr\} \,\mathrm{d} y +\vk{P}\bigl\{\bigl\llvert Z_i- y^{1/\gamma}
\bigr\rrvert \le x_i, \ZW{i\in I} \bigr\} \bigr]\,\mathrm{d} y
\\[-2pt]
&& \quad = \int_{-x_k^{\gamma}}^{x_k^{\gamma}} \vk{P}\bigl\{\bigl\llvert
Z_i+ \sign (y) \llvert y \rrvert^{1/\gamma} \bigr\rrvert \le
x_i, \ZW{i\in I} \bigr\} \,\mathrm{d} y,
\end{eqnarray*}
hence the proof follows.
\end{pf*}

\begin{pf*}{Proof of Theorem~\ref{T2}} First, we show that $\X_n= \X_n^{(1)}, n\ge1$ is the
$k$-dimensional marginal of some $(k+1)$-dimensional\vadjust{\goodbreak}
elliptical random vector. Define therefore a new random vector $\vk
{Y}_n, n\ge1$ with stochastic representation
\[
\vk{Y}_n\stackrel{d} {=}R_{k+1} A_n^*
\vk{U}_{k+1},
\]
where $\vk{U}_{k+1}$ is uniformly distributed on $\kal{S}_{k+1}$
independent of $R_{k+1}\sim H_{k+1}$, and $A_n^*$ is a non-singular
$(k+1)$-dimensional square matrix. Choose $A_n^*,n\ge1$ such that $\SI^*_n = A_n^*(A_n^*)^\top$ is again a correlation
matrix satisfying
\[
\bigl(\SI_n^*\bigr)_{II} = \ZW{\SI_n}, \qquad
I = \{1 ,\ldots,k\}, J = \{ k+1\},
\]
and
\[
\lim_{n \to\infty} a_n^2 \bigl(\vk{1}
\vk{1}^\top- \SI_n^*\bigr)= \Gamma^*\in(0,\IF
)^{(k+1)\times(k+1)}, \qquad\bigl(\Gamma^*\bigr)_{II}= \Gamma, \vk{1}\in
\R^{k+1}.
\]
Since $\SI_n,\SI_n^*$ are positive definite, by condition \eqref
{eq:reis} this construction is possible. Note that $\SI_n^*$ satisfies
\eqref{eq:reis} with $c_n = 2 b_n/a_n$ and limit matrix $\Gamma^*\in
[0,\IF)^{(k+1)\times(k+1)}$. We write for notational simplicity
$ (\Gamma^*)_{IJ} = \boldsymbol{\theta}/2$ and assume that
$\boldsymbol{\theta}$
has positive components.
It is well-known (see Cambanis~\cite{Cambanis1981}) that
\[
\vk{U}_{k+1} \stackrel{d} {=}\bigl(\vk{U} W, \sqrt{1- W^2}
\kal{J}\bigr),
\]
with $W$ a positive random variable such that $W^2 \stackrel{d}{=}\kal
{B}_{k/2,1/2}$, and $\kal{J}$ a Bernoulli random variable taking
values $-1,1$
with equal to probability $1/2$. Furthermore $\kal{J}, \vk{U}, $ and
$W$ are mutually independent.

By the assumption, $R^2_k \stackrel{d}{=}(R_{k+1})^2 \kal{B}_{k/2,1/2}$
with $R_{k+1}\sim H_{k+1}$ independent of $\kal{B}_{k/2,1/2}$,
implying $\Y_{n,I} \stackrel{d}{=}\X_n.$ Since the distribution function
of $\X_n$ depends on $\Sigma_n$ and not on $A_n$, and further $\Sigma_n$ is positive definite we can assume that
$A_n$ is a lower
triangular matrix. We construct $A_n^*$ to be also a non-singular lower
triangular matrix. With the same notation as in the proof of
{Theorem~\ref{th:kabl:0}}
we have
%
\begin{eqnarray}
\label{alves} (\Y_{n})_I | Y_{n,k+1}=q_n(y)
&\stackrel{d} {=}& R_{y,n,k} B_n \vk{U} + \bigl(
\Sigma^*_{n}\bigr)_{IJ}q_n(y), \qquad n\ge1,
\end{eqnarray}
where $B_n$ is a lower triangular matrix satisfying $B_n B_n^\top=
\Sigma_n- (\Sigma_{n}^*)_{IJ}(\Sigma_{n}^*)_{JI} $,
and $R_{y,n,k}$, \mbox{$n\ge1$} (being independent of $\vk{U}$) has survival
function $\overline{Q}_{y,n,k+1}$ given by \eqref{qN}.
As in the proof of {Theorem~\ref{th:kabl:0}}
\[
R_{y,n,k} \stackrel{d} {\to}\kal{R}_k \sim
\kal{H}_{k}, \qquad n\to \IF.
\]
By \eqref{eq:reis} and the fact that $B_nB_n^\top$ (and not the
matrix $B_n$) defines the conditional distribution we can choose $B_n$
such that $\lim_{n \to\infty} a_n B_n= B$ with
$ BB^\top= (\vk{1} \boldsymbol{\theta}^\top+ \boldsymbol{\theta
}\vk{1}^\top-
\Gamma)/2. $ Hence, for any $\vk{x}\in(0,\IF)^k$ utilising further
\eqref{alves} we obtain
\begin{eqnarray*}
&& \lim_{n \to\infty}\vk{P}\{a_n \vk{L}_n >
\vk{x}\}
\\
&& \quad = \lim_{n \to\infty}\vk{P}\{\forall i \le k\dvt  a_n L_{ni}
> x_i \}
\\
&& \quad = \lim_{n \to\infty}\vk{P}\bigl\{\forall i \le k\dvt  a_n\llvert
X_{ni} \rrvert > x_i \bigr\}^n
\\
&& \quad = \exp \Bigl(- \lim_{n \to\infty} n \vk{P}\bigl\{\exists i \le k\dvt
a_n\llvert X_{ni} \rrvert \le x_i \bigr\}
\Bigr)
\\
&& \quad = \exp \biggl(- \lim_{n \to\infty} n \biggl[ \int_{0}^\IF
\vk {P}\bigl\{\exists i \le k\dvt  a_n\llvert Y_{ni} \rrvert
\le x_i| Y_{n,k+1}= y/a_n \bigr\} \,\mathrm{d}
G_n( y)
\\
&&\hphantom{\exp \biggl(- \lim_{n \to\infty} n \biggl[} \qquad {} + \int_{-\IF}^0 \vk{P}\bigl\{\exists i \le k\dvt
a_n\llvert Y_{ni} \rrvert \le x_i |
Y_{n,k+1}= y/a_n \bigr\} \,\mathrm{d} G_n( y) \biggr]
\biggr)
\\
&& \quad = \exp \biggl(- \lim_{n \to\infty} n \int_{0}^\IF
\bigl[ \vk{P}\bigl\{\exists i \le k\dvt  a_n\llvert Y_{ni}
\rrvert \le x_i | Y_{n,k+1}= y/a_n \bigr\}
\\
&& \hphantom{\exp \biggl(- \lim_{n \to\infty} n \int_{0}^\IF
\bigl[}\qquad {} + \vk{P}\bigl\{\exists i \le k\dvt  a_n\llvert Y_{ni}
\rrvert \le x_i | Y_{n,k+1}= -y/a_n \bigr\}
\bigr]\,\mathrm{d} G_n( y) \biggr)
\\
&& \quad = \exp \biggl(- \lim_{n \to\infty} n \int_{0}^\IF
\bigl[ \vk {P}\bigl\{\exists i \le k\dvt  a_n\llvert Z_{ni}+
d_{ni} y/a_n \rrvert \le x_i \bigr\} \\
&& \hphantom{\exp \biggl(- \lim_{n \to\infty} n \int_{0}^\IF
\bigl[}\qquad {}+ \vk{P}
\bigl\{\exists i \le k\dvt  a_n\llvert Z_{ni}-
d_{ni} y/a_n \rrvert \le x_i \bigr\} \bigr]\,\mathrm{d} G_n( y) \biggr)
\\
&& \quad = \exp \biggl(- \int_{0}^\IF \bigl[ \vk{P}
\bigl\{\exists i \le k\dvt  \llvert Z_{i}+ y \rrvert \le x_i
\bigr\} + \vk{P}\bigl\{\exists i \le k\dvt  \llvert Z_{i} - y \rrvert
\le x_i \bigr\} \bigr]\,\mathrm{d} y^{ \gamma} \biggr)
\\
&& \quad = \exp \biggl(- \int_{0}^\IF \bigl[ \vk{P}
\bigl\{\exists i \le k\dvt  \bigl\llvert Z_{i}+ y^{1/\gamma} \bigr
\rrvert \le x_i \bigr\} + \vk{P}\bigl\{\exists i \le k\dvt  \bigl\llvert
Z_{i} - y^{1/\gamma} \bigr\rrvert \le x_i \bigr\}
\bigr]\,\mathrm{d} y \biggr)
\\
&& \quad = \exp \biggl(- \int_{\R} \vk{P}\bigl\{\exists i \le k\dvt
\bigl\llvert Z_{i}+ \sign(y) \llvert y \rrvert^{1/\gamma} \bigr
\rrvert \le x_i \bigr\} \,\mathrm{d} y \biggr),
\end{eqnarray*}
with $(Z_1 ,\ldots,Z_k)^\top=\kal{R}_k B \vk{U}$, and thus the claim
follows.
\end{pf*}

\begin{pf*}{Proof of Theorem~\ref{theo:main}} By Theorem~4.1 in Hashorva
and Pakes \cite
{Hashorva2010}, $H \in\mathrm{GMDA} (w)$ is equivalent with $G\in\mathrm{GMDA} (w)$. Let $B_n,\vk{Y}_n,
n\ge1$ be as in the proof of {Theorem~\ref{T2}} and adopt below the same
notation as therein.
Conditioning on $Y_{n,k+1}=q_n(y)=a_ny+ b_n, $ with $y\in\R$ such that
$G(q_n(y))\in(0,1),n \ge1$ we have that
\eqref{alves} holds, with $R_{y,n,k} $ independent of $\vk{U}$
satisfying (see Hashorva~\cite{Hashorva2009a})
\[
\frac{1}{\sqrt{ a_n b_n}} R_{y,n,k} \stackrel{d} {\to}\kal{R}, \qquad n\to\IF,
\]
where $\kal{R}^2\sim\chi^2_{k+1}$, and $\kal{R}_k>0$.
Next, $G \in\mathrm{GMDA} (w)$, \eqref{eq:reis} and the choice of $B_n$
imply for any $\vk{x}\in\R^k $ (omitting some details)
\begin{eqnarray*}
 &&\lim_{n \to\infty}\vk{P}\{\vk{M}_{n}\le a_n
\vk{x} +b_n \vk{1}\}
\\
&& \quad = \lim_{n \to\infty}\bigl[ 1- \vk{P}\bigl\{\exists i \le k\dvt  X_{ni}>
q_n(x_i) \bigr\} \bigr]^n
\\
&& \quad = \exp \Bigl(- \lim_{n \to\infty} n \vk{P}\bigl\{\exists i \le k\dvt
X_{ni} > q_n(x_i) \bigr\} \Bigr)
\\
&& \quad = \exp \biggl(- \lim_{n \to\infty} n\int_{\R} \vk{P}\bigl
\{\exists i \le k\dvt  Y_{ni} > q_n(x_i)|
Y_{n,k+1}= q_n(y) \bigr\} \,\mathrm{d} G\bigl(q_n(y)
\bigr) \biggr)
\\
&& \quad = \exp \biggl(- \lim_{n \to\infty} n\int_{\R} \vk{P}
\biggl\{\exists i \le k\dvt  \frac{1}{\sqrt{ a_n b_n}} R_{y,n,k} \bigl([
\sqrt{b_n/a_n}B_n] \vk{U}\bigr)_i
\\
&& \hphantom{\exp \biggl(- \lim_{n \to\infty} n\int_{\R} \vk{P}
\biggl\{}\qquad {}> x_i-y d_{ni} + [1- d_{ni}]b_n/a_n
\biggr\}\,\mathrm{d} G\bigl(q_n(y)\bigr) \biggr)
\\
&& \quad = \exp \biggl(- \int_{\R} \vk{P}\{\exists i \le k\dvt
Z_i > x_i-y + \theta_i/2 \} \exp(-y)\,\mathrm{d}y
\biggr),
\end{eqnarray*}
with $\vk{Z} \eqsim \mathit{Gauss}[\Gamma]$. Recall $\kal{R}_k \vk{U}$ is a
$k$-dimensional Gaussian random vector with independent components, and
further note that
the choice of $\theta_i$ above is arbitrary. The assumption that
\eqref{it} holds also for $m=k$ needed to define $\Y_n$ can now be
dropped since the limit
distribution is independent of that assumption, and further the
convergence in distribution holds without imposing that assumption,
hence the proof
is complete.
\end{pf*}

\begin{pf*}{Proof of Theorem~\ref{theo:Weib}} First note that
Theorem~4.5 in Hashorva~\cite{Hashorva2010} states that $H \in\mathrm{WMDA} (\alpha), \alpha>0$ is equivalent with $G\in\mathrm{WMDA}
(\alpha+ (k-1)/2)$. We proceed
as in the proof of {Theorem~\ref{theo:main}} (keeping the same notation).
Conditioning on the event $Y_{n,k+1}=q_n(y)= 1- a_ny, $ with $y$ such
that $G(q_n(y))\in(0,1),n \ge1$ and constants $a_n$ defined in
\eqref{an}
we have that again \eqref{alves} holds. In view of Hashorva \cite
{Hashorva2009a} for any $y>0$ 
\[
\frac{1}{\sqrt{ a_n }} R_{y,n,k} \stackrel{d} {\to}\sqrt{2y} \widetilde{\kal
{R}_\alpha}, \qquad n\to\IF,
\]
with $\widetilde{\kal{R}_\alpha}\sim\widetilde{\kal{H}}_\alpha$
where $\widetilde{\kal{H}}_\alpha(0)=0$ and $\widetilde{\kal
{R}_\alpha}^2 \stackrel{d}{=}
\kal{B}_{k/2,\alpha}$.
Furthermore
\[
\lim_{u \to\infty} \frac{1-G(1- x/u)}{1- G(1-1/u)}= x^{\alpha+
(k-1)/2}\qquad \forall x\in(0,
\IF)
\]
holds. Hence for any $\vk{x}\in(-\IF,0)^k $, we obtain (set $G_n(y)=
G(1- a_n y)$)
\begin{eqnarray*}
 &&\lim_{n \to\infty}\vk{P}\{\vk{M}_{n} \le
\vk{1}+a_n\vk {x}\}
\\
&& \quad = \exp \biggl(- \lim_{n \to\infty} n\int_{0}^\IF
\vk{P}\biggl\{\exists i \le k\dvt  \frac{1}{\sqrt{a_n}} R_{y,n,k} \biggl(
\frac{B_n}{\sqrt{a_n}} \vk {U}\biggr)_i \\
&& \hphantom{\exp \biggl(- \lim_{n \to\infty} n\int_{0}^\IF
\vk{P}\biggl\{}\qquad {}> x_i+y
d_{ni}+ [1- d_{ni}]/a_n \biggr\}\,\mathrm{d}
G_n(y) \biggr)
\\
&& \quad = \exp \biggl(- \int_{0}^\IF\vk{P}\bigl\{
\exists i \le k\dvt  Z_i > [ x_i+y + \theta_i/2
]/\sqrt{2 y}\bigr\} \,\mathrm{d} y^{\alpha+ (k-1)/2} \biggr),
\end{eqnarray*}
with $\vk{Z}\eqsim\mathfrak{E}[\Gamma;\widetilde{\kal{H}}_\alpha
]$, and thus the proof is complete.\vadjust{\goodbreak}
\end{pf*}

\begin{pf*}{Proof of Theorem~\ref{th:last}}
\begin{longlist}[(B)]
\item[(A)] Let $G$ denote the
distribution function of
$S_1 Y_{11}(1)$, and let $\Phi$ denote the standard Gaussian
distribution function on $\R$. The Mills ratio
asymptotics (see, e.g., Lu and Li~\cite{Lu2009}) implies $Y_{11}(1)\in
\kal{W}( 1/\sqrt{2
\uppi}, -1, 1/2,2)$. Consequently, by Lemma~2.1 in Arendarczyk and D\c
{e}bicki~\cite{Arendarczyk2011}
\begin{eqnarray*}
1- G(x)&=& \bigl(1+\mathrm{o}(1)\bigr) \biggl(\frac{2\uppi}{2+p_1} \biggr)%
^{1/2}
\frac{C_1}{\sqrt{2 \uppi}} A^{-\alpha_{1}}x^{\fracf{\alpha
_{1}(p_1-1)+p_{1}}{2+p_{1}}} \\
&&{}\times\exp \bigl( -
\bigl(L_{1}A^{-p_{1}}+A^{2}/2\bigr)x^{\fracc{2p_{1}}{%
2+p_{1}}}
\bigr) 
\\
&=& \bigl(1+\mathrm{o}(1)\bigr)\frac{C_1}{\sqrt{2+p_1}} A^{-\alpha_{1}}x^{\fracf{\alpha
_{1}(p_1-1)+p_{1}}{2+p_{1}}} \exp
\bigl(B x^{\fracc{2p_{1}}{%
2+p_{1}}}\bigr)
, \qquad x\to\IF,
\end{eqnarray*}
with $A = (p_1L_1)^{1/(2+p_1)}, B = L_{1}A^{-p_{1}}+A^{2}/2>0 .$
Hence $G\in\mathrm{GMDA} (w)$ with
\[
w(x)= B \frac{2p_1}{2+ p_1} x^{(p_1- 2)/(2+p_1)}, \qquad x>0.
\]
Set $b_n = G^{-1}(1- 1/n), n>1$ with $G^{-1}$ the generalised inverse
of $G$. Now, by \eqref{Davis}
%
\begin{eqnarray}
\label{Davis2} \lim_{n \to\infty} \frac{b_n}{b_n^*}&=&1,
\end{eqnarray}
where $b_n^* = \Psi^{-1}(1-1/n),n>1$ and $\Psi$ is some distribution
function satisfying
\[
1- \Psi(x) = \bigl(1+\mathrm{o}(1)\bigr) \exp\bigl(-B x^{\fracc{2p_{1}}{2+p_{1}}}\bigr), \qquad x\to
\IF.
\]
The above asymptotics implies
%
\begin{eqnarray}
\label{eq:gumb2} \lim_{n \to\infty} n \bigl( 1- G(a_nx+
b_n) \bigr)&=& \exp(- x) \qquad \forall x\in\R,
\end{eqnarray}
with
\[
b_n=\bigl(1+\mathrm{o}(1)\bigr) \biggl( \frac{\ln n}{B}
\biggr)^{(2+p_1)/(2 p_1)} , \qquad a_n= \frac
{1}{w(b_n)}=
\frac{ (2+p_1)b_n^{(2-p_1)/(2+p_1)}}{2 p_1 B}, \qquad n\to\IF.
\]
Consequently, as $n\to\IF$
\[
\frac{b_n}{a_n} = \bigl(1+\mathrm{o}(1)\bigr) \frac{2 p_1}{2+ p_1} \ln n ,
\]
hence \eqref{claim:M} follows by Theorem~3.1 of Kabluchko \cite
{Kabluchko2011} and {Theorem~\ref{theo:main}}.

\item[(B)] Since $\Phi\in\mathrm{GMDA}(w)$ with scaling function $w(x)=x,x>0$
Theorem~3 in Hashorva~\cite{Hashorva2009b} implies
\begin{eqnarray*}
1- G(x)&=& \bigl(1+\mathrm{o}(1)\bigr) \Gamma(\alpha+1) \vk{P}\bigl\{S> 1- 1/\bigl(x w(x)
\bigr)\bigr\} \vk {P}\bigl\{Y_{11}(1)> x\bigr\} , \qquad x\to\IF
\end{eqnarray*}
and thus $G\in\mathrm{GMDA} (w)$. If $a_n,b_n, n\ge1$ are defined by
\eqref{eq:gumb2}, then
Theorem~3.1 in Kabluchko~\cite{Kabluchko2011} and {Theorem~\ref{theo:main}}
establishes \eqref{claim:M}. By the form of $w(\cdot)$ we have
$\lim_{n \to\infty} a_nb_n=1$, and further
\eqref{Davis2}
holds with $b_n^* = \Phi^{-1}(1-1/n),n>1$. Consequently,
$b_n=(1+\mathrm{o}(1))\sqrt{2 \ln n}$ for all large $n$, and thus the result
follows.\hfill\qed
\end{longlist}\noqed
\end{pf*}

\begin{pf*}{Proof of Theorem~\ref{th:last:2}} Let $\vk{S}_n^{(i)}$ and
$\X_n^{(i)}, i \le n,
n\ge1$ be such that
$S_{nj}^{(i)}=S_{ni}(t_j), t_j\in\R, j\le k $ and $\X_{n}^{(i)},
i\le
n$ are independent copies of the Gaussian random vector
$X_{n1}(t_j),1 \le j\le k$.
By the assumptions of the theorem, the proof follows if we show that
the limit of the minima of the absolute values for the
triangular array
$ \vk{S}_n^{(i)} \X_n^{(i)}, i \le n, n\ge1$ converges to the random
vector $\boldsymbol{\kal{L}}$ such that
\begin{eqnarray*}
\vk{P}\{\boldsymbol{\kal{L}} > \vk{x}\} &=& \exp \biggl(- \int_{\R}
\vk {P}\bigl\{\exists i \le k\dvt  S_i\llvert y+Z_i \rrvert
\le x_i \bigr\} \,\mathrm{d} y \biggr), \qquad\vk{x}\in(0,\IF)^k,
\end{eqnarray*}
where $\vk{S}:=\vk{S}_{1}^{(1)}$ is independent centered Gaussian
random vector $\vk{Z}$ with incremental variance matrix $\Gamma$
which has components $\gamma_{ij}=\Gamma(t_i,t_j)$. The proof follows
with similar arguments as that of {Theorem~\ref{T2}} since $\vk{S}_n^{(i)}$
is, by the assumption, independent of $\X_n^{(i)}$.
\end{pf*}

\section*{Acknowledgements}
I would like to thank Zakhar Kabluchko, the referees of the paper, an
Associate Editor and the Editor-in-Chief, Professor Richard Davis
for several valuable comments and various corrections which improved
both the mathematics and the presentation.
Support by the Swiss National Science Foundation Grant 200021-134785 is
kindly acknowledged.


%

\printhistory


\begin{thebibliography}{29}

\bibitem{Arendarczyk2011}
%
\begin{barticle}[mr]
\bauthor{\bsnm{Arendarczyk},~\bfnm{Marek}\binits{M.}} \AND
\bauthor{\bsnm{D{\c{e}}bicki},~\bfnm{Krzysztof}\binits{K.}}
(\byear{2011}).
\btitle{Asymptotics of supremum distribution of a {G}aussian process
over a
{W}eibullian time}.
\bjournal{Bernoulli}
\bvolume{17}
\bpages{194--210}.
\bid{doi={10.3150/10-BEJ266}, issn={1350-7265}, mr={2797988}}
\bptok{imsref}%
\end{barticle}
%
\endbibitem

\bibitem{Bogachev1998}
%
\begin{bbook}[mr]
\bauthor{\bsnm{Bogachev},~\bfnm{Vladimir~I.}\binits{V.I.}}
(\byear{1998}).
\btitle{Gaussian Measures}.
\bseries{Mathematical Surveys and Monographs}
\bvolume{62}.
\baddress{Providence, RI}: \bpublisher{Amer. Math. Soc.}
\bid{mr={1642391}}
\bptok{imsref}%
\end{bbook}
%
\endbibitem

\bibitem{Breiman1965}
%
\begin{barticle}[mr]
\bauthor{\bsnm{Breiman},~\bfnm{L.}\binits{L.}}
(\byear{1965}).
\btitle{On some limit theorems similar to the arc-sin law}.
\bjournal{Theory   Probab. Appl.}
\bvolume{10}
\bpages{323--331}.
\bptok{imsref}%
\end{barticle}
%
\endbibitem

\bibitem{Brown1977}
%
\begin{barticle}[mr]
\bauthor{\bsnm{Brown},~\bfnm{Bruce~M.}\binits{B.M.}} \AND
\bauthor{\bsnm{Resnick},~\bfnm{Sidney~I.}\binits{S.I.}}
(\byear{1977}).
\btitle{Extreme values of independent stochastic processes}.
\bjournal{J. Appl. Probability}
\bvolume{14}
\bpages{732--739}.
\bid{issn={0021-9002}, mr={0517438}}
\bptok{imsref}%
\end{barticle}
%
\endbibitem

\bibitem{Cambanis1981}
%
\begin{barticle}[mr]
\bauthor{\bsnm{Cambanis},~\bfnm{Stamatis}\binits{S.}},
\bauthor{\bsnm{Huang},~\bfnm{Steel}\binits{S.}} \AND
\bauthor{\bsnm{Simons},~\bfnm{Gordon}\binits{G.}}
(\byear{1981}).
\btitle{On the theory of elliptically contoured distributions}.
\bjournal{J. Multivariate Anal.}
\bvolume{11}
\bpages{368--385}.
\bid{doi={10.1016/0047-259X(81)90082-8}, issn={0047-259X}, mr={0629795}}
\bptok{imsref}%
\end{barticle}
%
\endbibitem

\bibitem{Davis1988}
%
\begin{barticle}[mr]
\bauthor{\bsnm{Davis},~\bfnm{Richard}\binits{R.}} \AND
\bauthor{\bsnm{Resnick},~\bfnm{Sidney}\binits{S.}}
(\byear{1988}).
\btitle{Extremes of moving averages of random variables from the
domain of
attraction of the double exponential distribution}.
\bjournal{Stochastic Process. Appl.}
\bvolume{30}
\bpages{41--68}.
\bid{doi={10.1016/0304-4149(88)90075-0}, issn={0304-4149}, mr={0968165}}
\bptok{imsref}%
\end{barticle}
%
\endbibitem

\bibitem{Davis2008}
%
\begin{barticle}[mr]
\bauthor{\bsnm{Davis},~\bfnm{Richard~A.}\binits{R.A.}} \AND
\bauthor{\bsnm{Mikosch},~\bfnm{Thomas}\binits{T.}}
(\byear{2008}).
\btitle{Extreme value theory for space-time processes with heavy-tailed
distributions}.
\bjournal{Stochastic Process. Appl.}
\bvolume{118}
\bpages{560--584}.
\bid{doi={10.1016/j.spa.2007.06.001}, issn={0304-4149}, mr={2394763}}
\bptok{imsref}%
\end{barticle}
%
\endbibitem

\bibitem{De2006}
%
\begin{bbook}[mr]
\bauthor{\bparticle{de} \bsnm{Haan},~\bfnm{Laurens}\binits{L.}}
\AND
\bauthor{\bsnm{Ferreira},~\bfnm{Ana}\binits{A.}}
(\byear{2006}).
\btitle{Extreme Value Theory: An Introduction}.
\bseries{Springer Series in Operations Research and Financial Engineering}.
\baddress{New York}: \bpublisher{Springer}.
\bid{mr={2234156}}
\bptok{imsref}%
\end{bbook}
%
\endbibitem

\bibitem{Embrechts1997}
%
\begin{bbook}[mr]
\bauthor{\bsnm{Embrechts},~\bfnm{Paul}\binits{P.}},
\bauthor{\bsnm{Kl{\"u}ppelberg},~\bfnm{Claudia}\binits{C.}} \AND
\bauthor{\bsnm{Mikosch},~\bfnm{Thomas}\binits{T.}}
(\byear{1997}).
\btitle{Modelling Extremal Events: For Insurance and Finance}.
\bseries{Applications of Mathematics (New York)}
\bvolume{33}.
\baddress{Berlin}: \bpublisher{Springer}.
\bid{mr={1458613}}
\bptok{imsref}%
\end{bbook}
%
\endbibitem

\bibitem{Falk2010}
%
\begin{bbook}[auto:STB|2012/08/29|15:03:31]
\bauthor{\bsnm{Falk},~\bfnm{M.}\binits{M.}},
\bauthor{\bsnm{H{\"u}sler},~\bfnm{J.}\binits{J.}} \AND
\bauthor{\bsnm{Reiss},~\bfnm{R.~D.}\binits{R.D.}}
(\byear{2010}).
\btitle{Laws of Small Numbers: Extremes and Rare Events}, \bedition
{3rd} ed.
\baddress{Basel}: \bpublisher{Birkh\"auser}.
\bptok{imsref}%
\end{bbook}
%
\endbibitem



\bibitem{Hashorva2005}
%
\begin{barticle}[mr]
\bauthor{\bsnm{Hashorva},~\bfnm{Enkelejd}\binits{E.}}
(\byear{2005}).
\btitle{On the max-domain of attractions of bivariate elliptical arrays}.
\bjournal{Extremes}
\bvolume{8}
\bpages{225--233}.
\bid{doi={10.1007/s10687-006-7969-6}, issn={1386-1999}, mr={2275920}}
\bptok{imsref}%
\end{barticle}
%
\endbibitem

\bibitem{Hashorva2006}
%
\begin{barticle}[mr]
\bauthor{\bsnm{Hashorva},~\bfnm{Enkelejd}\binits{E.}}
(\byear{2006}).
\btitle{On the multivariate {H}\"usler--{R}eiss distribution attracting the
maxima of elliptical triangular arrays}.
\bjournal{Statist. Probab. Lett.}
\bvolume{76}
\bpages{2027--2035}.
\bid{doi={10.1016/j.spl.2006.05.022}, issn={0167-7152}, mr={2329248}}
\bptok{imsref}%
\end{barticle}
%
\endbibitem

\bibitem{Hashorva2008}
%
\begin{barticle}[mr]
\bauthor{\bsnm{Hashorva},~\bfnm{Enkelejd}\binits{E.}}
(\byear{2008}).
\btitle{On the max-domain of attraction of type-{III} elliptical triangular
arrays}.
\bjournal{Comm. Statist. Theory Methods}
\bvolume{37}
\bpages{1543--1551}.
\bid{doi={10.1080/03610920801893871}, issn={0361-0926}, mr={2440452}}
\bptok{imsref}%
\end{barticle}
%
\endbibitem

\bibitem{Hashorva2009b}
%
\begin{barticle}[mr]
\bauthor{\bsnm{Hashorva},~\bfnm{Enkelejd}\binits{E.}}
(\byear{2009}).
\btitle{Conditional limit results for type {I} polar distributions}.
\bjournal{Extremes}
\bvolume{12}
\bpages{239--263}.
\bid{doi={10.1007/s10687-008-0078-y}, issn={1386-1999}, mr={2533952}}
\bptok{imsref}%
\end{barticle}
%
\endbibitem

\bibitem{Hashorva2009a}
%
\begin{barticle}[mr]
\bauthor{\bsnm{Hashorva},~\bfnm{Enkelejd}\binits{E.}}
(\byear{2009}).
\btitle{Conditional limits of {$W\sb p$}-scale mixture distributions}.
\bjournal{J. Statist. Plann. Inference}
\bvolume{139}
\bpages{3501--3511}.
\bid{doi={10.1016/j.jspi.2009.04.001}, issn={0378-3758}, mr={2549098}}
\bptok{imsref}%
\end{barticle}
%
\endbibitem

\bibitem{Hashorva2011}
%
\begin{barticle}[auto:STB|2012/08/29|15:03:31]
\bauthor{\bsnm{Hashorva},~\bfnm{E.}\binits{E.}}
(\byear{2013}).
\btitle{On Beta-product convolutions}.
\bjournal{Scand. Actuar. J.}
\bvolume{2013}
\bpages{69--83}.
\bptok{imsref}%
\end{barticle}
%
\endbibitem

\bibitem{Hashorva2010}
%
\begin{barticle}[mr]
\bauthor{\bsnm{Hashorva},~\bfnm{Enkelejd}\binits{E.}} \AND
\bauthor{\bsnm{Pakes},~\bfnm{Anthony~G.}\binits{A.G.}}
(\byear{2010}).
\btitle{Tail asymptotics under beta random scaling}.
\bjournal{J.~Math. Anal. Appl.}
\bvolume{372}
\bpages{496--514}.
\bid{doi={10.1016/j.jmaa.2010.07.045}, issn={0022-247X}, mr={2678878}}
\bptok{imsref}%
\end{barticle}
%
\endbibitem

\bibitem{H2011a}
%
\begin{barticle}[mr]
\bauthor{\bsnm{H{\"u}sler},~\bfnm{J{\"u}rg}\binits{J.}},
\bauthor{\bsnm{Piterbarg},~\bfnm{Vladimir}\binits{V.}} \AND
\bauthor{\bsnm{Rumyantseva},~\bfnm{Ekaterina}\binits{E.}}
(\byear{2011}).
\btitle{Extremes of {G}aussian processes with a smooth random variance}.
\bjournal{Stochastic Process. Appl.}
\bvolume{121}
\bpages{2592--2605}.
\bid{doi={10.1016/j.spa.2011.06.006}, issn={0304-4149}, mr={2832415}}
\bptok{imsref}%
\end{barticle}
%
\endbibitem

\bibitem{H2011b}
%
\begin{barticle}[mr]
\bauthor{\bsnm{H{\"u}sler},~\bfnm{J{\"u}rg}\binits{J.}},
\bauthor{\bsnm{Piterbarg},~\bfnm{Vladimir}\binits{V.}} \AND
\bauthor{\bsnm{Zhang},~\bfnm{Yueming}\binits{Y.}}
(\byear{2011}).
\btitle{Extremes of {G}aussian processes with random variance}.
\bjournal{Electron. J. Probab.}
\bvolume{16}
\bpages{1254--1280}.
\bid{doi={10.1214/EJP.v16-904}, issn={1083-6489}, mr={2827458}}
\bptok{imsref}%
\end{barticle}
%
\endbibitem

\bibitem{H1989}
%
\begin{barticle}[mr]
\bauthor{\bsnm{H{\"u}sler},~\bfnm{J{\"u}rg}\binits{J.}} \AND
\bauthor{\bsnm{Reiss},~\bfnm{Rolf-Dieter}\binits{R.D.}}
(\byear{1989}).
\btitle{Maxima of normal random vectors: Between independence and complete
dependence}.
\bjournal{Statist. Probab. Lett.}
\bvolume{7}
\bpages{283--286}.
\bid{doi={10.1016/0167-7152(89)90106-5}, issn={0167-7152}, mr={0980699}}
\bptok{imsref}%
\end{barticle}
%
\endbibitem

\bibitem{Jessen2006}
%
\begin{barticle}[mr]
\bauthor{\bsnm{Jessen},~\bfnm{Anders~Hedegaard}\binits{A.H.}} \AND
\bauthor{\bsnm{Mikosch},~\bfnm{Thomas}\binits{T.}}
(\byear{2006}).
\btitle{Regularly varying functions}.
\bjournal{Publ. Inst. Math. (Beograd) (N.S.)}
\bvolume{80}
\bpages{171--192}.
\bid{doi={10.2298/PIM0694171J}, issn={0350-1302}, mr={2281913}}
\bptok{imsref}%
\end{barticle}
%
\endbibitem

\bibitem{Kabluchko2011}
%
\begin{barticle}[mr]
\bauthor{\bsnm{Kabluchko},~\bfnm{Zakhar}\binits{Z.}}
(\byear{2011}).
\btitle{Extremes of independent {G}aussian processes}.
\bjournal{Extremes}
\bvolume{14}
\bpages{285--310}.
\bid{doi={10.1007/s10687-010-0110-x}, issn={1386-1999}, mr={2824498}}
\bptok{imsref}%
\end{barticle}
%
\endbibitem

\bibitem{Kabluchko2009}
%
\begin{barticle}[mr]
\bauthor{\bsnm{Kabluchko},~\bfnm{Zakhar}\binits{Z.}},
\bauthor{\bsnm{Schlather},~\bfnm{Martin}\binits{M.}} \AND
\bauthor{\bparticle{de} \bsnm{Haan},~\bfnm{Laurens}\binits{L.}}
(\byear{2009}).
\btitle{Stationary max-stable fields associated to negative definite
functions}.
\bjournal{Ann. Probab.}
\bvolume{37}
\bpages{2042--2065}.
\bid{doi={10.1214/09-AOP455}, issn={0091-1798}, mr={2561440}}
\bptok{imsref}%
\end{barticle}
%
\endbibitem

\bibitem{Lu2009}
%
\begin{barticle}[mr]
\bauthor{\bsnm{Lu},~\bfnm{Dawei}\binits{D.}} \AND
\bauthor{\bsnm{Li},~\bfnm{Wenbo~V.}\binits{W.V.}}
(\byear{2009}).
\btitle{A note on multivariate {G}aussian estimates}.
\bjournal{J. Math. Anal. Appl.}
\bvolume{354}
\bpages{704--707}.
\bid{doi={10.1016/j.jmaa.2009.01.046}, issn={0022-247X}, mr={2515252}}
\bptok{imsref}%
\end{barticle}
%
\endbibitem

\bibitem{Omey2010}
%
\begin{bincollection}[mr]
\bauthor{\bsnm{Omey},~\bfnm{Edward}\binits{E.}} \AND
\bauthor{\bsnm{Segers},~\bfnm{Johan}\binits{J.}}
(\byear{2010}).
\btitle{Generalised regular variation of arbitrary order}.
In \bbooktitle{Stability in Probability}.
\bseries{Banach Center Publ.}
\bvolume{90}
\bpages{111--137}.
\baddress{Warsaw}: \bpublisher{Polish Acad. Sci. Inst. Math.}
\bid{doi={10.4064/bc90-0-8}, mr={2798855}}
\bptok{imsref}%
\end{bincollection}
%
\endbibitem

\bibitem{Penrose1991}
%
\begin{barticle}[mr]
\bauthor{\bsnm{Penrose},~\bfnm{Mathew~D.}\binits{M.D.}}
(\byear{1991}).
\btitle{Minima of independent {B}essel processes and of distances between
{B}rownian particles}.
\bjournal{J. London Math. Soc. (2)}
\bvolume{43}
\bpages{355--366}.
\bid{doi={10.1112/jlms/s2-43.2.355}, issn={0024-6107}, mr={1111592}}
\bptok{imsref}%
\end{barticle}
%
\endbibitem


\bibitem{Sibuya1960}
%
\begin{barticle}[mr]
\bauthor{\bsnm{Sibuya},~\bfnm{Masaaki}\binits{M.}}
(\byear{1960}).
\btitle{Bivariate extreme statistics. {I}}.
\bjournal{Ann. Inst. Statist. Math. Tokyo}
\bvolume{11}
\bpages{195--210}.
\bid{issn={0020-3157}, mr={0115241}}
\bptok{imsref}%
\end{barticle}
%
\endbibitem

\end{thebibliography}
\end{document}